\documentclass[11pt]{article}
\usepackage{subfigure} 
\usepackage{epsfig}
\usepackage{amsmath, amsthm}

\newcommand{\be}{\begin{eqnarray}}     	\newcommand{\ee}{\end{eqnarray}}

\newcommand{\vol}{\mathrm{Vol}}
\newcommand{\dist}{\mathrm{dist}}
\newcommand{\rem}{\mathrm{Rm}}
\newcommand{\ric}{\mathrm{Ric}}
\newcommand{\diam}{\mathrm{diam}}

\title{Properties of the solutions of the conjugate heat equations}

\author{Richard Hamilton, Natasa Sesum}

\date{} 
\theoremstyle{plain}
\newtheorem{dummy}{Dummy}

\theoremstyle{definition}
\newtheorem{definition}[dummy]{Definition}

\theoremstyle{plain}

\newtheorem{lemma}[dummy]{Lemma}
\newtheorem{theorem}[dummy]{Theorem}
\newtheorem{proposition}[dummy]{Proposition}
	
\newtheorem{claim}[dummy]{Claim}

\begin{document}

\maketitle

\begin{abstract}
In this paper we consider the class $\mathcal{A}$ of those solutions
$u(x,t)$ to the conjugate heat equation $\frac{d}{dt}u = -\Delta u +
Ru$ on compact K\"ahler manifolds $M$ with $c_1 > 0$ (where $g(t)$ changes
by the unnormalized K\"ahler Ricci flow, blowing up at $T <
\infty$), which satisfy Perelman's differential Harnack inequality
on $[0,T)$. We show $\mathcal{A}$ is nonempty. If $|\ric(g(t))| \le
\frac{C}{T-t}$, which is alaways true if we have type I singularity,
we prove the solution $u(x,t)$ satisfies the elliptic type Harnack
inequlity, with the constants that are uniform in time.  If the flow
$g(t)$ has a type I singularity at $T$, then $\mathcal{A}$ has
excatly one element.
\end{abstract}

\begin{section}{Introduction}

Let $M$ be a K\"ahler manifold with $c_1(M) > 0$, of complex
dimension $n$. Consider the solutions to the unnormalized K\"ahler
Ricci flow,
\begin{equation}
\label{equation-flow}
\frac{d}{dt}g_{i\bar{j}} = -R_{i\bar{j}}.
\end{equation}
It is known in the case of the unnormalized K\"ahler Ricci flow that
it shrinks  to a point, after some finite time $T < \infty$.
Let $T' < T$ and let $u = (4\pi(T'-t))^{-n}e^{-f}$ satisfy the conjugate heat equation
\begin{equation}
\label{equation-conjugate}
\frac{d}{dt}u = -\Delta u + Ru.
\end{equation}
This implies $f$ satisfies,
\begin{equation}
\label{equation-evolution-f}
\frac{d}{dt}f = -\Delta f + |\nabla f|^2 - R + \frac{n}{T-t}.
\end{equation}
Let 
\begin{equation}
\label{equation-v}
v = [(T'-t)(2\Delta f - |\nabla f|^2 + R) + f - 2n]u,
\end{equation}
which is such that $\int_M v$ is exactly Perelman's functional
$\mathcal{W}$. He proved it is monotonically increasing along the flow,
that is,
$$\frac{d}{dt}\mathcal{W} = 2\tau\int_M |R_{i\bar{j}} + \nabla_i\nabla_{\bar{j}}f - g_{i\bar{j}}|^2udV \ge 0,$$ 
If $u$ tends to a $\delta$-function as $t\to T'$, in
\cite{perelman2002}, Perelman proved $v \le 0$ for all $t\in [0,T']$.
He also proved that under the same assumptions as above, for any
smooth curve $\gamma(t)$ in $M$,
\begin{equation}
\label{equation-T'}
-\frac{d}{dt}f(\gamma(t),t) \le \frac{1}{2}(R(\gamma(t),t) + 
|\dot{\gamma}(t)|^2) - \frac{1}{2(T'-t)}f(\gamma(t),t),
\end{equation}
for all $t\in [0,T']$.

\begin{definition}
\label{definition-admissable}
We will say that a smooth function $f$ is {\it admissable} if for
any smooth curve $\gamma(t)$ in $M$, 
\begin{equation}
\label{equation-diff-harnack}
-\frac{d}{dt}f(\gamma(t),t) \le \frac{1}{2}(R(\gamma(t),t) + 
|\dot{\gamma}(t)|^2) - \frac{1}{2(T-t)}f(\gamma(t),t),
\end{equation}
for all $t\in [0,T)$, while the metric changes by the Ricci flow equation 
(\ref{equation-flow}) and $T$ is a time at which the flow disappears.
\end{definition}

We will prove the following results about $u$.

\begin{theorem}
\label{theorem-harnack}
If $|\ric(g(t))| \le \frac{c}{T-t}$, which translates to the
condition, $\ric \ge -c$ along the normalized K\"ahler Ricci flow, the
set $\mathcal{A}$ is nonempty and there is a uniform constant $C$, so
that,
$$\max_{M\times[0,T)}u(x,t) \le C\min_{M\times[0,T)}u(x,t).$$
If we assume the flow has a type I singularity, meaning that $|\rem(g(t))| \le \frac{C}{T-t}$,
there is excatly one element in $\mathcal{A}$, that is, the solution
to the conjugate heat equation (\ref{equation-conjugate}), existing
all the way up to $T$ and satisfying (\ref{equation-diff-harnack})
is unique. 
\end{theorem} 

The organization of the paper is as follows. In section $2$ we will
give the proof of Theorem \ref{theorem-harnack}. Complex two
dimensional case will be discussed in section $3$. In section $4$ we
will discuss Perelman's reduced distance function and show how its
definition can be extended to a distance function with a base point at
$(p,T)$, where $p$ is the point to which the flow shrinks and $T$ is
the singular time.
\end{section}

\begin{section}{Harnack type estimates and the uniqueness of $u$}
\label{section-harnack}

We will assume for the moment that $\mathcal{A}$ is not empty and prove
that each element $u(t)\in \mathcal{A}$ satisfies the elliptic Harnack inequality
at each time slice with the uniform constant, not depending on time, and 
that such a solution is unique if $g(t)$ has type I singularity at $T$.

\begin{proposition}
\label{proposition-harnack}
If $|\ric(g(t))| \le \frac{C}{T-t}$ along the flow $g(t)$, there exists
a uniform constant $\tilde{C}$, so that 
$$\max_{M}u(x,t) \le \tilde{C}\min_M u(x,t),$$
for all $t\in [0,T)$.
\end{proposition}

\begin{proof}
Take $t_1 < t_2 < T$ and $x_1,x_2\in M$. Let $\gamma(t)$ be a curve
that will be chosen later, so that it connects $x_1$ and
$x_2$, that is, $\gamma(t_1) = x_1$ and $\gamma(t_2) = x_2$.
Since $f$ is an admissable function, it satisfies (\ref{equation-diff-harnack}),
for a chosen curve $\gamma$. Integrate it in $t\in [t_1,t_2]$. We get,
\begin{equation}
\label{equation-harnack-applied}
f(x_1,t_1)\sqrt{T-t_1} \le f(x_2,t_2)\sqrt{T-t_2} + 
\frac{1}{2}\int_{t_1}^{t_2}\sqrt{T-t}(R(\gamma(t),t) + |\dot{\gamma}(t)|^2)dt.
\end{equation}
By translation in time, we may assume $T = 1/2$. It easily follows
that if we rescale the flow, that is, if $\tilde{g}(s(t)) =
\frac{g(t)}{1-t/T} = Tg(t)/(T-t)$ with $s(t) = -T\ln (1 - t/T)$, we
get a normalized K\"ahler Ricci flow, satisfying,
$$\frac{d}{ds}\tilde{g} = \tilde{g} - \ric(\tilde{g}),$$
for all $s\in [0,\infty)$. By Perelman's results,
$|R(\tilde{g}(s))| \le C$ and $\diam(M,\tilde{g}(s)) \le C$ along the flow. 
This implies
\begin{equation}
\label{equation-scalar}
R(g(t)) \le \frac{C}{T-t},
\end{equation}
and
\begin{equation}
\label{equation-diameter}
\diam(M,g(t)) \le C\sqrt{T-t}.
\end{equation}
As a matter of reparametrization, we also get 
\begin{equation}
\label{equation-volume}
\vol_{g(t)} = (1-t/T)^n\vol_{\tilde{g}(s(t))}(M) = C(T-t)^n.
\end{equation}
We will estimate the integral term appearing in (\ref{equation-harnack-applied}). By
(\ref{equation-scalar}) we have,
\begin{equation}
\label{equation-est1}
\int_{t_1}^{t_2}\sqrt{T-t} R(\gamma(t),t) dt \le C\int_{t_1}^{t_2}\frac{dt}{\sqrt{T-t}}
= C\frac{t_2 - t_1}{\sqrt{T-t_1} + \sqrt{T-t_2}}.
\end{equation}
Without loosing a generality assume that $|\ric(g(t))|(T-t) \le 1$,
which may be always achieved by rescaling. We have the simple claim
that follows immediatelly from the evolution equation for $g(t)$.

\begin{claim}
If $\ric(g(t))(T-t) \ge -g(t)$ for all $t\in [0,T)$, for any $0 \le t < s < T$, we have,
$$g(s) \le \frac{T-t}{T-s}g(t).$$
In particular, for any vector $v$, we have
$$|v|^2_{g(s)} \le \frac{T-t}{T-s}|v|^2_{g(t)}.$$
\end{claim}

Let $s_1 = s(t_1)$.
If we choose $\gamma$ to be a minimal geodesic from
$x_1$ to $x_2$ with respect to $g(t_1)$, by (\ref{equation-diameter}) and by 
the previous claim, we have,
\begin{eqnarray}
\label{equation-est2}
\int_{t_1}^{t_2}\sqrt{T-t}|\dot{\gamma}|^2dt &=& 
\int_{T-t_2}^{T-t_1}\sqrt{\tau}|\dot{\gamma}|^2_{g(T-\tau)}d\tau  \\
&\stackrel{\sqrt{\tau}=s}{=}& \int_{\sqrt{T-t_2}}^{\sqrt{T-t_1}}\frac{1}{2}|\gamma'|^2_{g(T-s^2)}ds \nonumber \\
&\le& C\int_{\sqrt{T-t_2}}^{\sqrt{T-t_1}}\frac{T-t_1}{T-t_2}|\gamma'|^2_{g(t_1)}ds \nonumber \\
&=& \tilde{C}\frac{T-t_1}{T-t_2}\frac{\dist^2_{g(t_1)}(x_1,x_2)}{\sqrt{T-t_1} - \sqrt{T-t_2}} \nonumber \\
&\stackrel{(\ref{equation-diameter})}{\le}& \bar{C}\frac{(T-t_1)^2\sqrt{T-t_1}}{(T-t_2)(t_2-t_1)} \nonumber.
\end{eqnarray}
From estimates (\ref{equation-harnack-applied}), (\ref{equation-est1})
and (\ref{equation-est2}), we get
$$f(x_1,t_1) \le \frac{\sqrt{T-t_2}}{\sqrt{T-t_1}}f(x_2,t_2) + 
C\frac{t_2-t_1}{T-t_1} + C\frac{(T-t_1)^2}{(t_2-t_1)(T-t_2)}.$$
If (*) $t_2-t_1 \sim T-t_1$, e.g. $t_2-t_1 = T-t_2 = \delta$, then,
\begin{equation}
\label{equation-compare}
f(x_1,t_1) \le \frac{1}{\sqrt{2}}f(x_2,t_2) + C,
\end{equation}
and since $x_1, x_2$ were two arbitrary points, we have
\begin{equation}
\label{equation-min-max}
\max_M f(\cdot,t_1) \le \frac{1}{\sqrt{2}}\min_M f(\cdot,t_2) + C.
\end{equation}

We claim there is some $\tilde{A}$ so that $f(x,t) \ge -\tilde{A}$, for all $x\in M$ and all $t\in [0,T)$.
Assume $\max_M f(\cdot,t) \le -A$. That implies $f(x,t) \le -A$, for all $x\in M$ and
therefore,
$$(4\pi\tau(T-t))^n = \int_M e^{-f}dV_t \ge e^A\vol_t(M) = e^A C(T-t)^n,$$
for a uniform constant $C$, which is not possible for big enough $A$ (notice 
that the bigness of $A$ does not depend on $t\in [0,T)$). Estimate
(\ref{equation-min-max}) now implies
$$\min_Mf(\cdot,t_2) \ge -\tilde{A},$$
for a uniform constant $\tilde{A}$, independent of $t_2$. Rewrite 
(\ref{equation-compare}) as
\begin{eqnarray*}
\max_M f(\cdot,t_1) &\le& \frac{1}{\sqrt{2}}(\min_M f(\cdot,t_2) + \tilde{A}) - 
\frac{\tilde{A}}{\sqrt{2}} + C \\
&\le& \min_M f(\cdot,t_2) + \tilde{A} + C_1 = \min_M f(\cdot,t_2) + C_2.
\end{eqnarray*}
If we denote by $M(t) = \max_M u(\cdot,t)$ and by $m(t) = \min_M u(\cdot,t)$, this yields
\begin{equation}
\label{equation-first}
M(t_2) \le C(\frac{T-t_1}{T-t_2})^n  m(t_1)\le \tilde{C}m(t_1).
\end{equation}
The evolution equation $\frac{d}{dt}u = -\Delta u + Ru$ at the points where $u(\cdot,t)$
achieves its maximum becomes,
$$\frac{d}{dt}M(t) \ge RM(t) \ge -\frac{C}{T-t}M(t),$$
which yields $M(t_2) \ge (\frac{T-t_2}{T-t_1})^C M(t_1) = \tilde{C}M(t_1)$. 
This together with (\ref{equation-first}) gives,
\begin{equation}
\label{equation-final}
M(t_1) \le Cm(t_1),
\end{equation}
for a uniform constant $C$ and all $t_1\in [0,T)$, which is an analogue of 
the Harnack inequality that we have in the elliptic case.
\end{proof}

We will now prove the nonemptiness of $\mathcal{A}$.
 
\begin{lemma}
\label{lemma-nonempty}
For every K\"ahler manifold $M$ and every unnormalized K\"ahler Ricci
flow $g(t)$, $|\mathcal{A}| \ge 1$.
\end{lemma}

\begin{proof}
For every unnormalized K\"ahler Ricci flow $g(t)$ there is a finite
time $T$ at which the flow disappears. Take an arbitrary increasing
sequence of times $t_i \uparrow T$ and a sequence of points $p_i\in
M$. For every $i$, let $u_i(t) = (4\pi(t_i-t))^{-n}e^{-f_i(t)}$ be a
solution to the conjugate heat equation (\ref{equation-conjugate}),
such that $u_i(t)$ converges to a $\delta$-function as $t\to t_i$,
concentrated at $p_i$. Let $v_i(t)$ be a corresponding $v$-function
as in (\ref{equation-v}). Due to Perelman (see \cite{perelman2002}),
we have $v_i(t) \le 0$ for all $t\in [0,t_i]$ and for every smooth
curve $\gamma(t)$,
$$-\frac{d}{dt}f_i(\gamma(t),t) \le \frac{1}{2}(R(\gamma(t),t) + 
|\dot{\gamma}(t)|^2) - \frac{1}{2(t_i-t)}f_i(\gamma(t),t),$$
holds for all $t\in [0,t_i]$.

Fix $t_j$ from the sequence of times and consider $u_i(t)$ for $i\ge
j$.  Then the metrics $\{g(t)\}$ are uniformly equivalent and their
geometries are uniformly bounded for $t\in [0,t_j]$, with bounds
depending on $g(t)|_{[0,t_j]}$. Each $u_i(t)$ satisfies the conjugate
heat equation (\ref{equation-conjugate}), with $\int_M u_i(t) dV_t = 1$.
Moreover,
$$\frac{d}{dt}f_i = -\Delta f_i + |\nabla f_i|^2 - R + \frac{n}{t_i-t},$$  
which implies,
$$\frac{d}{dt}(f_i)_{\min} \le \frac{C}{T-t} + \frac{n}{t_i-t} \le
\frac{\tilde{C}}{t_i-t},$$
and therefore, for $t\in [0,t_j]$,
$$(f_i)_{\min}(t_j) \le (f_i)_{\min}(t) + \tilde{C}\ln(\frac{t_i-t}{t_i-t_j}).$$
This yields,
$$(f_i)_{\min}(t) \ge -C(t_j),$$
for $i$ big enough, so that $t_i \ge \frac{T+t_j}{2}$, and henceforth,
\begin{equation}
\label{equation-bound-u-i}
\max_{M\times[0,t_j]}u_i(x,t) \le e^{-C(t_j)}(2\pi(T-t_j))^{-n} = \tilde{C}(t_j).
\end{equation}
If we denote by $\tilde{u}_i(x,t) = u_i(x,t)^{1/2}$, integrating 
$v_i(x,t) \le 0$, similarly as in \cite{notes}, by using H\"older and Sobolev
inequalities, 
$$\int_M |\nabla \tilde{u}_i|^2 dV_t \le C_j,$$
which together with (\ref{equation-bound-u-i}) imply,
$$\sup_{t\in [0,t_j]}||u_i(t)||_{W^{1,2}} \le C_j,$$
for all $i$ big enough. By standard parabolic estimates applied to
$u_i(t)$ and $t\in [0,t_j]$, it follows there exists $C(k,l,n,t_j)$ so that
$$|\frac{d^l}{dt^l}u_i(t)|_{C^k} \le C(k,l,n,t_j),$$
for all $t\in [0,t_j]$ and all $i$ sufficiently big. Extract a subsequence
$u_i(t,x)$ that converges in $C^{k,l}(M\times[0,t_j])$ norm to some
function $u(t,x)$, defined on $[0,t_j]$ that continues to be a solution
to the conjugate heat equation (\ref{equation-conjugate}). By taking
larger and larger $j$, diagonalizing our sequence $u_i(t)$, taking into
account the uniqueness of the limit, we get a function $u(t,x)$,
defined on $M\times[0,T)$ and a subsequence $u_i(t,x)$, so that
$u_i(t,x)\stackrel{C^{k,l}(M\times[0,T'])}{\to} u(t,x)$, for every $T'<T$.
Moreover, 
\begin{itemize}
\item
$u(t,x)$ satisfies the conjugate heat equation (\ref{equation-conjugate}) for
all $t\in [0,T)$.
\item
$u(t,x) = (4\pi(T-t))^{-n}e^{-f}$, where $f$ is an admissable function in the
sense of Definition \ref{definition-admissable}.  
\end{itemize}
In particular, this implies $|\mathcal{A}| \ge 1$
\end{proof} 

We will now prove the uniqueness part of Theorem \ref{theorem-harnack}.

\begin{proposition}
\label{proposition-uniqueness}
If $g(t)$ has type I singularity at $T$, then $|\mathcal{A}| = 1$, that is, the 
solution $u$ to the conjugate heat equation $\frac{d}{dt}u = -\Delta u + Ru$,
that satisfies Perelman's differential Harnack inequality (\ref{equation-diff-harnack})
is unique.
\end{proposition}

\begin{proof}
Assume there are at least two different solutions $u_1$ and $u_2$,
with the properties as above. By Theorem \ref{theorem-harnack}, we have
$M_j(t) \le Cm_j(t)$, for all $t\in [0,T)$ and $j\in \{1,2\}$. We will omit
the subscript $j$ below. On the other hand, we have 
the integral normalization condition $(4\pi(T-t))^{-n}\int_M u dV_{g(t)} = 1$.
Combining these two facts, we get
\begin{equation}
\label{equation-low-u}
m(t) \ge \frac{1}{C\vol_{g(t)}(M)},
\end{equation}
\begin{equation}
\label{equation-up-u}
M(t) \le C\frac{1}{\vol_{g(t)}(M)},
\end{equation}
where a constant $C$ comes from (\ref{equation-final}).
Take a sequence $t_i\to T$ and consider a sequence of rescaled metrics
$g_i(t) = (T-t_i)^{-1}g((T-t_i)t + t_i)$, for $t\in [-(T-t_i)^{-1}t_i,1)$.
We get
\begin{eqnarray*}
|\rem(g_i(t))| &=& |\rem(g(t(T-t_i)+t_i)|(T-t_i) \\
&=& |\rem(g(t(T-t_i)+t_i)|\frac{(T-t_i)}{T-t_i-t(T-t_i)} \\
&\le& C\frac{1}{1-t} \le \tilde{C},
\end{eqnarray*}
for all $t\in [-(T-t_i)^{-1}t_i, 1/2]$.
By Perelman's volume noncollapsing result and by Hamilton's
compactness theorem, there is a subsequence $(M,g_i(t))$, 
converging to another solution to the ancient K\"ahler Ricci 
solution $(M,h(t))$, defined for $t\in (-\infty, 1/2]$. If we 
also rescale our solution $u$, together with our metric $g(t)$,
estimates (\ref{equation-low-u}), (\ref{equation-up-u})
and (\ref{equation-volume}) give
\begin{eqnarray*}
u_i(t) &=& (T-t_i)^n u(t(T-t_i)+t_i) \\
&\le& C\frac{(T-t_i)^n}{\vol_{g(t(T-t_i)+t_i)}(M)} \\
&=& \tilde{C}\frac{(T-t_i)^n}{(T-t_i)^n(1-t)^n} \le \bar{C},
\end{eqnarray*}
for all $t\in [-(T-t_i)^{-1}t_i,1/2]$. Notice that function $f$ rescales as
$f_i(t) = f(t_i+t(T-t_i))$. Similarly we get the uniform lower
bound on $u_i(t)$, that is, there exists a uniform constant $\tilde{C} > 1$
so that
\begin{equation}
\label{equation-est-uniform}
\frac{1}{\tilde{C}} \le u_i(t) \le \tilde{C},
\end{equation}
on $M\times [-(T-t_i)^{-1}t_i,1/2]$. 
Functions $u_i(t)$ satisfy backward parabolic equations
$$\frac{d}{dt}u_i(t) = -\Delta u_i(t) + R(g_i(t))u_i(t).$$
Uniform $C^0$ estimate (\ref{equation-est-uniform}) imply,
$$\int_{-1}^{1/2}\int_M|\nabla u_i|^2dV_{g_i(t)}dt \le C.$$
Since $|\nabla u_i| = (4\pi(1-t))^{-n}|\nabla f_i|u_i$ and 
(\ref{equation-est-uniform}), we have,
$$\int_{-1}^{1/2}\int_M|\nabla f_i|^2dV_{g_i(t)}dt \le \tilde{C}.$$
By standard parabolic estimates applied to $f_i$ satisfying 
(\ref{equation-evolution-f}), we have,  
\begin{equation}
\label{equation-est-f}
\sup_{t\in [-1,1/2]}|f_i(t)|_{C^k} \le C(k).
\end{equation}

\begin{claim}
\label{claim-boundness-W}
There is a uniform constant $C$ so that $\sup_{t\in
[0,T)}\mathcal{W}(g(t),f(t),T-t) \le C$.
\end{claim}

\begin{proof}
Consider our sequence $t_i\to T$. Then, by the scale invariance of $\mathcal{W}$
and by estimates (\ref{equation-est-f}), 
\begin{eqnarray}
\label{equation-upper-W}
\mathcal{W}(g(t_i),f(t_i),T-t_i) &=& \mathcal{W}(g_i(0),f_i(0),1) \nonumber \\
&=& (4\pi)^{-n}\int_M (R_i + |\nabla f_i|^2 + f_i - 2n)e^{-f_i}dV_{g_i(0)} \nonumber \\
&\le& C.
\end{eqnarray}
By Perelman's monotonicity formula, $\mathcal{W}(g(t),f(t),T-t)$ increases 
in time, which together with (\ref{equation-upper-W}) imply the statement of the 
claim.
\end{proof}
The previous claim and Perelman's monotonicity formula for $\mathcal{W}$ yield
the existence of a finite limit, $\lim_{t\to T} \mathcal{W}(g(t),f(t),T-t)$.
Let $a_i = \frac{T+t_i}{2}$. From before, we have that $g_i(t) \to h(t)$.
From our estimates (\ref{equation-est-f}) on $f_i(s)$, by extracting a subsequence
we may assume $f_i(s) \stackrel{C^k(M\times[-1,1/2]}{\to} f_h(s)$. We also have,
\begin{eqnarray}
\label{equation-important}
& &\mathcal{W}(g(a_i),f(a_i),T-a_i) - \mathcal{W}(g(t_i),f(t_i),T-t_i) =
\int_{t_i}^{a_i}\frac{d}{dt}\mathcal{W}dt = \\
&=& \int_{t_i}^{a_i}(4\pi(T-t))^{-n}\int_M (2(T-t)|R_{p\bar{q}} + 
\nabla_p\nabla_{\bar{q}}f - g_{p\bar{q}}|^2e^{-f}dV_{g(t)}dt + \nonumber \\
&+& |\nabla_p\nabla_q f|^2 + |\nabla_{\bar{p}}\nabla_{\bar{q}}f|^2)e^{-f}dV_{g(t)}dt \nonumber \\
&\ge& (4\pi(T-t_i)^{-n}\int_{t_i}^{a_i}\int_M ((T-t_i)(|R_{p\bar{q}} + 
\nabla_p\nabla_{\bar{q}}f - g_{p\bar{q}}|^2 +
|\nabla_p\nabla_q f|^2 + \nonumber \\
&+& |\nabla_{\bar{p}}\nabla_{\bar{q}}f|^2) e^{-f}dV_{g(t)}dt \nonumber \\
&=& (4\pi)^{-n}\int_0^{1/2}\int_M(|\ric_i + \nabla\bar{\nabla}f_i - g_i|^2 + |\nabla\nabla f_i|^2
+ |\bar{\nabla}\bar{\nabla}f_i|^2)e^{-f_i}dV_{g_i(s)}ds. \nonumber 
\end{eqnarray}
The left hand side of (\ref{equation-important}) converges to zero, while its right hand side
converges to 
$$\int_0^{1/2}\int_M (|\ric(h) + \nabla\bar{\nabla}f_h - h|^2 + |\nabla\nabla f_h|^2 +
|\bar{\nabla}\bar{\nabla}f_h|^2)e^{-f_h}dV_{h(s)}ds.$$
This yields $h(s)$ is a K\"ahler Ricci soliton and it satisfies,
\begin{eqnarray*}
R_{p\bar{q}}(h) + \nabla_p\nabla_{\bar{q}}f_h - h_{p\bar{q}} &=& 0, \\
f_{\bar{p}\bar{q}} = f_{pq} &=& 0.
\end{eqnarray*}
In other words, what we get is the following: if we have two different
solutions $u_1 = (4\pi(T-t))^{-n}e^{-f_1}$ and $u_2 =
(4\pi(T-t))^{-n}e^{-f_2}$, to each of them we can apply the reasoning
from above. We can consider $u^1_i(t) = (T-t_i)^{-n}u_1(t_i+t(T-t_i))$
and $u^2_i(t) = (T-t_i)^{-n}u_2(t_i+t(T-t_i))$ and as above, we can
conclude $f_1(t_i+t(T-t_i))\stackrel{C^k}{\to} f^1_h$ and
$f_2(t_i+t(T-t_i))\stackrel{C^k}{\to} f^2_h$, where $f^1_h$ 
and
$f^2_h$ both satisfy,
$$\ric(h) + \nabla\bar{\nabla}f^1_h - h = 0,$$
$$\ric(h) + \nabla\bar{\nabla}f^2_h - h = 0.$$
This implies $\Delta f_h^1 = \Delta f_h^2$, which yields $f_h^1 = f_h^2 + C$, for 
some constant $C$. Since $\int_M e^{-f_h^1}dV_h = \int_M e^{-f_h^2}dV_h$, we get $C=0$.
This means, for $t\in [-1,1/2]$, 
$$\frac{u_1(t_i+t(T-t_i))}{u_2(t_i+t(T-t_i))} = \frac{u_i^1(t)}{u_i^2(t)} 
\stackrel{C^k}{\to} 1,$$
and in particular, by putting $t=0$,
\begin{equation}
\label{equation-conv-1}
\frac{u_1(t_i)}{u_2(t_i)} = \frac{u_i^1(0)}{u_i^2(0)} 
\stackrel{C^k}{\to} 1,
\end{equation}
as $i\to\infty$, where $u_i^1(t) = (T-t_i)^nu_1(t_i+t(T-t_i))$,
and similarly for $u_i^2(t)$.

A simple computation shows that the evolution equation for
$\frac{u_1}{u_2}$, since both functions $u_1$ and $u_2$ satisfy the
conjugate heat equation (\ref{equation-conjugate}) is
\begin{equation}
\label{equation-maximum}
\frac{d}{dt}\frac{u_1}{u_2} = -\Delta \frac{u_1}{u_2} - 
\nabla\ln u_2\nabla(\frac{u_1}{u_2}).
\end{equation}
If there is a time $t_0\in [0,T)$ and $x\in M$ so that $\frac{u_1(x,t_0)}{u_2(x,t_0)} \ge 1 + \delta$,
with $\delta > 0$, then $\max_M(\frac{u_1(\cdot,t_0)}{u_2(\cdot,t_0)}) \ge 1 + \delta$.
By the maximum principle applied to (\ref{equation-maximum}), we get 
$\max_M(\frac{u_1(\cdot,t)}{u_2(\cdot,t)})$ increases in time and therefore,
$$\max_M(\frac{u_1(\cdot,t)}{u_2(\cdot,t)}) \ge 1+\delta,$$
for all $t\in [t_0,T)$. This contradicts (\ref{equation-conv-1}) and
henceforth $u_1(t) = u_2(t)$ for all $t\in [0,T)$.
\end{proof}

\begin{proof}[Proof of Theorem \ref{theorem-harnack}]
The proof of the theorem follows from Lemma \ref{lemma-nonempty}, Proposition
\ref{proposition-harnack} and Proposition \ref{proposition-uniqueness}.
\end{proof}

\end{section}

\begin{section}{More on the uniqueness of $u$ for $n=2$}

In this section we will consider two dimensional K\"ahler, compact
manifolds $M$, with $c_1 > 0$. Let $g(t)$ be the K\"ahler Ricci flow
on such a manifold. In the K\"ahler case, the curvature integral
$\int_M|\rem|^2dV$ is always bounded in terms of topological
invariants, the first and the second Chern class. This integral is
scale invariant for $n=2$ which implies its significant importance
in that case. For example, using that in \cite{natasa-kr}
the following result has been proved.

\begin{theorem}
\label{theorem-convergence-two}
Let $g(t)$ be the normalized K\"ahler Ricci flow on a manifold as
above, with uniformly bounded Ricci curvatures along the flow. Then
for every sequence $t_i\to\infty$, there is a subsequence, so that
$(M,g(t_i+t))\to (M_{\infty},g_{\infty}(t))$, where $M_{\infty}$ is
the orbifold with finitely many isolated singularities and
$g_{\infty}(t)$ is a singular metric that satisfies the K\"ahler Ricci
soliton equation outside those singular points.
\end{theorem}

Combining Theorem \ref{theorem-convergence-two} and Theorem
\ref{theorem-harnack} yields the following result.

\begin{theorem}
\label{theorem-uniqueness-two}
If $g(t)$ is the unnormalized K\"ahler Ricci flow on a manifold $M$
as above, such that $|\ric(g(t))| \le \frac{C}{T-t}$, for a uniform 
constant $C$, there is a unique solution of the conjugate heat equation
(\ref{equation-conjugate}).
\end{theorem}

\begin{proof}
The proof is analogous to the proof of Theorem \ref{theorem-harnack},
since the only singularities we get in  two dimensional case are just
isolated points. Adopt the notation from the proof of Theorem
\ref{theorem-harnack}. For the rescaled sequence of metrics
$g_i(t) = (T-t_i)^{-1}g(t_i + t(T-t_i))$, we have that
$$\sup_{M\times[-t_i(T-t_i)^{-1},1)}|\ric(g_i(t))| \le C,$$
$$\delta < u_i(t) \le C,$$
for all $t\in [-t_i(T-t_i)^{-1},\frac{1}{2}]$ and therefore,
$$\int_{-1}{1/2}\int_M|\nabla f_i(s)|^2dV_{g_i(s)}ds \le C.$$
The last estimate implies that for every $i$, there is
$s_i\in [-1,1/2]$, so that
$$\int_M|\nabla f_i(s_i)|^2dV_{g_i(s_i)} \le C.$$
In the proof of Claim \ref{claim-boundness-W}, to prove the boundness
of $\mathcal{W}$, instead of considering
$\mathcal{W}(g(t_i),f(t_i),T-t_i)$ we will consider
$\mathcal{W}(g(t_i+s_i(T-t_i)),f(t_i+s_i(T-t_i)),(T-t_i)(1-s_i))$, for
$s_i\in [-1,1/2]$. The rest of the proof is same.  We also have the
same estimate (\ref{equation-important}) as before, where the left
hand side tends to zero as $i\to\infty$, due to the monotonicity and
the boundness of $\mathcal{W}$. Assume $p_1,\dots, p_N$ are the
singular points we get by taking the limit of the sequence
$(M,g(t_i+t))$, and that $p_1^i,\dots,p_N^i$ are the curvature
concentration points that are responsible for obtaining our
singularities in the limit. Let $\{D_j\}$ be the compact
exhaustion of $M_{\infty}\backslash\{p_1,\dots,p_N\}$. Our geometries
$g(t)$ are uniformly bounded on each of $D_j$, (those bounds
deteriorate when $j\to\infty$, that is, when we approach singularities).
Henceforth, the estimate (\ref{equation-important}) tells us we can extract
a subsequence, such that $(M,g(t_i+t))\to (M_{\infty},h(t))$ and 
$h(t)$ satisfies the K\"ahler Ricci soliton equation,
\begin{equation}
\label{equation-soliton-0}
\ric(h(t)) + \nabla\bar{\nabla}f_h(t) - h(t) = 0,
\end{equation}
away from singular points. As in Proposition
\ref{proposition-uniqueness}, if we assume there are at least two
different solutions of the conjugate heat equation, we will get at
least two different functions $f_h(t)$ and $f'_{h}(t)$, that satisfy
(\ref{equation-soliton-0}) away from singular points. Without loss of
generality assume there is only one singular point $p$.

\begin{claim}
Functions $f_h(t)$ and $f'_{h}t)$ coincide on $M_{\infty}\backslash\{p\}$
\end{claim}

\begin{proof}
Choose a sequence $\eta_k\to 1$ on $M_{\infty}$ with $\int|\nabla\eta_k|^2\to 0$ 
as $k\to\infty$, e.g.,
$$\eta_k(t) = \left \{
\begin{array}{lll}
0, & \mbox{for } x\in B(p,1/k^2), \\
1 - \frac{\ln(k^2\dist_h(p,x))}{\ln k}, & \mbox{for } x\in B(p,1/k)\backslash B(p,1/k^2), \\
1, & \mbox{for } x\in M_{\infty}\backslash B(p,1/k).
\end{array}
\right .$$ 
Denote by $F = f_h - f_h'$. It satisfies, $\Delta_h F = 0$ away from $p$.
Multiply it by $F^2\eta_k^2$ and then integrate over $M$. We get,
\begin{eqnarray*}
\int |\nabla F|^2\eta_k^2 dV_h&=& -\int \nabla\eta_k\eta_kF\nabla FdV_h \\
&\le& 4\int |F|^2|\nabla\eta_k|^2dV_h + \frac{1}{4}\int \eta_k^2|\nabla F|^2dV_h, \\
&\le& C\int |\nabla\eta_k|^2dV_h + \frac{1}{4}\int \eta_k^2|\nabla F|^2dV_h 
\end{eqnarray*}
which after taking $k\to\infty$ implies,
$$\int_{M_{\infty}\backslash\{p\}}|\nabla F|^2dV_h = 0.$$

As in \cite{anderson1989}, \cite{bando1989} and \cite{tian1990} one can show
$M_{\infty}\backslash \{p\}$ is connected. This amounts to having $F = const$ on 
$M_{\infty}\backslash \{p\}$. Becuase of the integral normalization condition for
$f_h$ and $f_h'$, we have $C=0$ and therefore, $f_h = f'_h$. 
\end{proof}

The rest of the proof
of Theorem \ref{theorem-uniqueness-two} is as in the proof of Theorem
\ref{theorem-harnack}.
\end{proof}
\end{section}

\begin{section}{Reduced distance function with the base point $(p,T)$}

In \cite{perelman2002} Perelman has introduced the reduced distance function
for the Ricci flow $g(t)$ defined for $t \in [0,T']$,
with respect to the base point $(p,T')$, for some $p\in M$, as follows.
For any point $q\in M$ and any $t\in [0,T']$, let
\begin{equation}
\label{equation-red-dist}
l(q,T'-t) = \frac{1}{2\sqrt{T'-t}}\int_t^{T'}\sqrt{T'-u}(R(\gamma(u),u) + 
|\dot{\gamma}|^2)du,
\end{equation}
where $\gamma$ is the $\mathcal{L}$-geodesic (minimizing the
integral in (\ref{equation-red-dist})), such that $\gamma(t) = q$ and
$\gamma(T') = p$. If in (\ref{equation-T'}) we choose $\gamma$ to be the
$\mathcal{L}$ geodesic, integrating (\ref{equation-T'}) in $t$, if
$\lim_{s\to T'}\sqrt{T'-s}f(p,s)$ (which is true if $u = 
(4\pi(T'-t))^{-n}e^{-f}$ tends to a $\delta$-function concentrated at $p$,
as $t\to T'$) yields, 
\begin{equation}
\label{equation-f-l}
f(q,t) \le l(q,T'-t).
\end{equation}

We would like to define some notion of the reduced distance for the 
K\"ahler Ricci flow, defined with respect to the base point $(p,T)$, where 
$p$ is a point to which our flow shrinks at singular time $T$.
The idea is roughly as follows.
Let $t_i\uparrow T$ as $i\to\infty$. Fix points $x,q\in M$ and $t\in [0,T)$; and for
each $i$ such that $t < t_i$, define the $\mathcal{L}$-distance with the base point
$(x,t_i)$, from $(x,t_i)$ to $(q,t)$ (as Perelman did in \cite{perelman2002}). We will
denote it by $L_i^x(q,t)$.
Define
$$\tilde{L}_i(q,t) = \inf_{x\in M}L_i^x(q,t),$$
where the infimum is taken over all base points $(x,t_i)$.
Assume $t_i \le t_j$. Take any base point $(x,t_j)$ and
let $\gamma_1$ be a restriction of $\gamma$ to time interval
$[t,t_i]$. Then,
\begin{eqnarray*}
L_j^x(q,t) &=& \int_t^{t_j}\sqrt{t_j-u}(R+|\dot{\gamma}|^2)du \ge
\int_t^{t_i}\sqrt{t_i-t}(R+|\dot{\gamma_1}|^2)du \\
&\ge& \tilde{L}_i(q,t),
\end{eqnarray*}
If we take the infimum over all base points $(x,t_j)$ in the previous 
inequality, we get
\begin{equation}
\label{equation-increasing}
\tilde{L}_j(q,t) \ge \tilde{L}_i(q,t),
\end{equation}
which means $\tilde{L}_i$ is an increasing sequence.
$$\tilde{L}_i(t,q) \le L_i^q(t,q) \le \int_t^{t_i}\sqrt{t_i-u}(R + |\dot{\gamma}|^2)du,$$
where we can take  $\gamma(t)$ to be a constant curve $\gamma(t) = q$.
Then $\dot{\gamma} = 0$. Due to Perelman, we have that for the K\"ahler Ricci flow
$\sup_{M\times[0,T)}|R|(T-t) \le C$ and therefore,
\begin{eqnarray}
\label{equation-bound}
\tilde{L}_i(t,q) &\le& \int_t^{t_i}\sqrt{t_i-u}Rdu \nonumber \\
&\le& \int_t^{t_i}\sqrt{t_i-u}\frac{C}{T-u}du \le 
\int_t^{t_i}\frac{1}{\sqrt{t_i-u}}du  \nonumber\\
&=& C\sqrt{t_i-t}.
\end{eqnarray}
By (\ref{equation-increasing}) and (\ref{equation-bound}) we get there
is a $\lim_{i\to\infty}\tilde{L}_i(q,t) = \tilde{L}(q,t)$. This
implies $\tilde{l}_i = \frac{1}{\sqrt{t_i-t}}\tilde{L}_i \to 
\frac{1}{\sqrt{T-t}}\tilde{L} = \tilde{l}$.
An estimate (\ref{equation-bound}) implies
$$\sup_{M\times[0,T)}\tilde{l}(q,t) \le C,$$
for a uniform constant $C$.

One interesting question would be whether $\tilde{l}_i$ satisfy
the similar inequalities to those that are satisfied by each of
$l_i^x$. Recall that Perelman has proved $l_i^x$ satisfies,
$$-(l_i^x)_t - \Delta l_i^x + |\nabla l_i^x|^2 - R + \frac{n}{t_i-t} \ge 
0,$$
$$2\Delta l_i^x - |\nabla l_i^x|^2 + R + \frac{l_i^x-2n}{t_i-t} \le 0.$$

{\bf Question:} Do the above inequalities persist after taking 
the infimum of $l_i^x$ over all $x\in M$?

If the answer to the above question were positive, this would yield
the monotonicity formula for $\tilde{V}(q) = (T-t)^{-n}\int_M 
e^{-\tilde{l}}dV_g$, for all $t\in [0,T)$.

\end{section}

\end{document}